\documentclass[prb,aps,epsf,eqsecnum]{revtex4}
\usepackage{graphicx}

\newcommand{\dx}{{{\mit \Delta} x}}
\newcommand{\dt}{{{\mit \Delta} t}}
\newcommand{\bu}{{\bf u}}
\newcommand{\lm}{\lambda}
\newcommand{\ld}{{\it r}}

\newcommand{\de}{\epsilon}
\newcommand{\e}{{\rm e}}
\newcommand{\ep}{{{\mit\Delta} t}}

\newcommand{\A}{{\bf A}}
\newcommand{\B}{{\bf B}}
\newcommand{\C}{{\bf C}}

\newcommand{\bea}{\begin{eqnarray}}
\newcommand{\eea}{\end{eqnarray}}
\newcommand{\baa}{\begin{array}}
\newcommand{\eaa}{\end{array}}
\newcommand{\be}{\begin{equation}}
\newcommand{\ee}{\end{equation}}
\newcommand{\ba}{\begin{eqnarray}}
\newcommand{\ea}{\end{eqnarray}}

\newcommand{\nn}{\nonumber}
\newcommand{\la}{\label}

\def\t1{e_{_T}}
\def\v1{e_{_V}}

\begin{document}
\title{Symplectic finite-difference methods for solving partial differential equations}

\author{Siu A. Chin}
\affiliation{Department of Physics, 
Texas A\&M University, College Station, TX 77843, USA}

\begin{abstract}
The usual explicit finite-difference method of solving partial differential equations
is limited in stability because it approximates the exact amplification factor by power-series.
By adapting the same exponential-splitting method of deriving symplectic integrators,
explicit symplectic finite-difference methods produce Saul'yev-type schemes which approximate
the exact amplification factor by rational-functions. As with conventional symplectic integrators, 
these symplectic finite-difference algorithms preserve 
important qualitative features of the exact solution. Thus the symplectic diffusing algorithm is 
{\it unconditionally stable} and the symplectic advection algorithm is {\it unitary}.
There is a one-to-one correspondence between symplectic integrators and symplectic
finite-difference methods, including the key idea that one can systematically improve an
algorithm by matching its modified Hamiltonian more closely to the original Hamiltonian.  
Consequently, the entire arsenal of symplectic 
integrators can be used to produce arbitrary high order time-marching algorithms for solving 
the diffusion and the advection equation.    

\end{abstract}
\maketitle

\section {Introduction}
The 1D	diffusion equation
\be
\frac{\partial u}{\partial t}=D\frac{\partial^2 u}{\partial x^2}
\la{deq}
\ee
can be solved numerically by applying the forward-time and central-difference 
approximations to yield the explicit algorithm
\be
u_j^\prime=u_j+r[ (u_{j+1}-u_j)-(u_j-u_{j-1})],
\la{fdone}
\ee
where $x_j=j\dx$, $u_j=u(x_j,t)$, $u_j^\prime=u(x_j,t+\dt)$ and
\be
\ld=\frac{\dt\,D}{\dx^2}.					
\ee
Under this (Euler) algorithm, each Fourier component $\tilde u_k=\e^{ikx}$ with wave number 
$k$ is amplified by a factor of
\be
g=1-4\ld\sin^2(k\dx/2),
\la{euler}
\ee
restricting stability ($|g|\le 1$) to 
the Courant-Friedrichs-Lewy\cite{cfl} (CFL) limit,
\be
\ld\le\frac12.
\ee 
Since explicit finite-difference methods approximate the exact amplification factor by 
power-series such as (\ref{euler}), it seems inevitable that they 
will eventually blow-up and be limited in stability. However, Saul'yev\cite{saul57,saul64} showed in the 50's
that, by simply replacing in (\ref{fdone}), either
\be  
(u_{j}-u_{j-1})\rightarrow (u_{j}^\prime-u_{j-1}^\prime)
\quad{\rm or}\quad
(u_{j+1}-u_j)\rightarrow (u_{j+1}^\prime-u_j^\prime)
\ee
one would have {\it unconditionally stable} algorithms:
\be
u_j^\prime=\beta\!_S\, u_{j-1}^\prime
+\gamma\!_S\, u_j+\beta\!_S\, u_{j+1},
\la{saula}
\ee
or 
\be
u_j^\prime=\beta\!_S\, u_{j-1}
+\gamma\!_S\, u_j+\beta\!_S\, u^\prime_{j+1},
\la{saulb}
\ee
where $\gamma\!_S$ and $\beta\!_S$ are Saul'yev's coefficients given by
\be
 \gamma\!_S=\frac{1-r}{1+r}\quad{\rm and}\quad \beta\!_S=\frac{r}{1+r}.
\ee
Algorithm (\ref{saula}) is explicit if it is evaluated in ascending order
in $j$ from left to right and if the left-most $u_1$ is a boundary value fixed in time. 
Similarly, algorithm (\ref{saulb}) is explicit if it is evaluated in descending order
in $j$ from right to left and if the right-most $u_N$ is a boundary value fixed in time.
Saul'yev also realized that both algorithms have large errors (including phase errors due to
their asymmetric forms), but if they are applied 
{\it alternately}, the error would be greatly reduced after such a pair-wise application.
This then gives rise to {\it alternating direction explicit} algorithms
advocated by Larkin\cite{lark64} and generalized to {\it alternating group explicit} 
algorithms by Evans\cite{evans83,evans85}.

There are four unanswered questions about Saul'yev asymmetric algorithms:
1) While it is easy to show that algorithm (\ref{saula}) and (\ref{saulb}) are 
unconditionally stable, there is no deeper understanding of this stability. 
2) The algorithms are not explicit in the case of periodic
boundary. What would be the algorithm if there are no fixed boundary values?
3) The alternating application of (\ref{saula}) and (\ref{saulb}) greatly reduces the
resulting error. How can one characterize this improvement precisely? 4) How can Saul'yev-type algorithms
be generalized to higher orders?

This work presents a new way of deriving finite-difference schemes based on exponential-splittings 
rather than Taylor expansions. Exponential-splitting
is the basis for developing symplectic integrators\cite{cre89,fr90,suzu90,yos90,yos93,hairer02,mcl02}, 
the hallmark of structure-preserving algorithms.
The finite-difference method presented here is properly ``symplectic" in the original sense that it has
certain ``intertwining" quality, resembling Hamilton's equations. It is also symplectic in the wider sense of 
structure-preserving, in that there is a Hamiltonian-like quantity that the algorithms seek to preserve. 

As will be shown, there is a one-to-one correspondence between
symplectic finite-difference methods and symplectic integrators. It is therefore useful
to summarize some basic results of symplectic integrators for later reference.  
Symplectic integrators are based on approximating
${\rm e}^{\de(\A+\B)}$ to any order in $\de$ via
a single product decomposition
 \be
{\rm e}^{\de( \A+ \B)}=
\prod_{i}
{\rm e}^{a_i\de \A}{\rm e}^{b_i\de \B},
\label{prod} 
\ee
where \A\ and \B\ are non-commuting operators (or matrices). Usually, \A+\B={\bf H} is
the Hamiltonian operator and ${\rm e}^{\de{\bf H}}$ is the evolution operator that evolves 
the system forward for time $\de$.
The key idea is to preserve the {\it exponential} character of the evolution operator.
The two first-order, Trotter\cite{trot58} approximations are
\be
T_{\rm 1A}(\de)={\rm e}^{\de \A}{\rm e}^{\de \B},\quad 
T_{\rm 1B}(\de)={\rm e}^{\de \B}{\rm e}^{\de \A},
\la{sym1ab}
\ee
and the two second-order Strang\cite{strang68} approximations are 
\ba
T_{\rm 2A}(\de)&=&T_{\rm 1A}(\de/2)T_{\rm 1B}(\de/2)
={\rm e}^{\frac12 \de \A}{\rm e}^{\de \B}{\rm e}^{\frac12 \de \A},\nn\\
T_{\rm 2B}(\de)&=&T_{\rm 1B}(\de/2)T_{\rm 1A}(\de/2)
={\rm e}^{\frac12 \de \B}{\rm e}^{\de \A}{\rm e}^{\frac12 \de \B}.
\la{sym2ab}
\ea
The approximation
\be
T_{\rm 2C}(\de)=\frac12\Bigl[T_{\rm 1A}(\de)+T_{\rm 1B}(\de)\Bigr]
\ee
is also second-order, but since it is no longer a single product of exponentials,
it is no longer symplectic.	In most cases, it is inferior to $T_{\rm 2A}$ and $T_{\rm 2B}$
because the time steps used in evaluating $T_{\rm 1A}$ and $T_{\rm 1B}$ are twice as large as those 
used in $T_{\rm 2A}$ and $T_{\rm 2B}$. 

Let $T_2$ denotes either $T_{\rm 2A}$ or $T_{\rm 2B}$. $T_2$ must be second order because for a 
left-right symmetric product as above, it must obey 
\be
T_2(-\de)T_2(\de)=1,
\la{two1}
\ee
and therefore must be of the form,  
\be
T_2(\de)={\rm e}^{\de {\bf H}+\de^3{\bf E}_3+\de^5{\bf E}_5+\cdots}
\la{oddstr}
\ee 
with only odd powers of $\de$ in the exponent. (Since there is no way for the operators in (\ref{two1}) 
to cancel if there are any even power terms in $\de$.) 
In (\ref{oddstr}), ${\bf E}_n$ denote higher order commutators of \A\ and {\B}.
The algorithm corresponding to $T_2$ then yields exact trajectories of  
the second-order {\it modified Hamiltonian}
\be
{\bf H}_2(\de)={\bf H}+\de^2{\bf E}_3+\de^4{\bf E}_5+\cdots .
\ee
A standard way of improving the efficiency of symplectic integrators is to generate a 
2n$^{th}$-order algorithm via a product of second-order 
algorithms\cite{cre89,suzu90,yos90}, via
\be
T_{2n}(\de)=\prod_{i=1}^{N} T_2(a_i\de).
\la{t2n}
\ee
Since the error structure of $T_2(\de)$ is given by (\ref{oddstr}), to preserve the original
Hamiltonian, one must choose $a_i$ to perserve the first power of $\de$,
\be
\sum_{i=1}^Na_i=1.
\la{sum1}
\ee
To obtain a fourth-order algorithm, one must eliminate the error term proportinal to $\de^3$ by
requiring, 
\be
\sum_{i=1}^Na_i^3=0.
\la{sum3}
\ee
For a sixth-order algorithm, one must require the above and
\be
\sum_{i=1}^Na_i^5=0,
\la{sum5}
\ee
and so on. While proofs of these assertions in terms of operators are not difficult, we will not need them.  
Symplectic finite-difference methods use a much simpler version of these ideas. Instead of
dealing with the evolution operator $\e^{\de{\bf H}}$, the finite difference method has
a proxy, the amplification factor, which is just a function. Order-conditions such as (\ref{sum3}) and
(\ref{sum5}) will then be obvious. Other results will be cited as needed, but these basic findings 
are sufficient to answer the four questions about Saul'yev's schemes. For the next two sections
we will give a detailed derivation of the symplectic diffusion and advection algorithms,
followed by a discussion of the diffusion-advection equation and a concluding summary.

\section {Symplectic diffusion algorithm}

Consider solving the diffusion equation (\ref{deq}) with periodic boundary condition
$u_{N+1}=u_1$ in the semi-discretized form,
\be
\frac{du_j}{dt}=\frac{D}{\dx^2}(u_{j+1}-2u_j+u_{j-1}).
\ee
Regarding $u_j$ as a vector, this is
\be
\frac{d\bu}{dt}=\A\bu,
\la{veq}
\ee
with
\be
\bu=\left(\begin{array}{c}
       u_1\\
       u_2\\
	   \vdots\\
	   u_N
      \end{array}\right), \qquad
\A=\frac{D}{\dx^2}
\left(\baa{ccccc}
-2 &  1  &   & &  1   \\
 1 &  -2 & 1 & &      \\
   & & \ddots &  &     \\
   &     &   1& -2& 1    \\
 1  &     &   & 1& -2    \\
\eaa\right),
\ee
and exact solution 
\be
\bu(t+\dt)=e^{\dt\A}\bu(t).
\la{expa}
\ee
The Euler algorithm corresponds to 
expanding out the exponential to first order in $\dt$
\be
\bu(t+\dt)=(1+\dt\A )\bu(t),
\la{fd1}
\ee
resulting in a power-series amplification factor (\ref{euler}), with limited stability.

If the exponential in (\ref{expa}) can be solved exactly, 
the amplification factor would be 
\be
g_{ex}=\e^{-h_{ex}},
\la{gex}
\ee
where
\be
h_{ex}=\ld4\sin^2(\theta/2) \quad{\rm and }\quad \theta\equiv k\dx.
\la{hex}
\ee 
The {\it amplification exponent} $h_{ex}$ here plays the role of a time
parameter $r$ times the original ``Hamiltonian" $h_0=4\sin^2(\theta/2)$.
The resulting algorithm will then be 
{\it unconditionally stable} for all $r>0$. In the limit of $\dx\rightarrow 0$,
each $k$-Fourier components will be damped by
$g_{ex}=\e^{-\dt D k^2}$, which is the exact solution to (\ref{deq}).

To preserve this important feature of the exact solution,  
one must seek alternative ways of approximating of $\e^{\dt \A}$ without doing any Taylor expansion. 
The structure of
$\A$ immediately suggests that it should decompose as
\be
\A=\sum_{j=1}^N\A_j,
\ee
where each $\A_j$ has only a single, non-vanishing $2\times 2$ matrix along the diagonal
connecting the $j$ and the $j+1$ elements:
\be
\A_j=\frac{D}{\dx^2}\left(\baa{cccccc}
  \ddots&    &   & &  &   \\
  & -1 & 1 & &  &    \\
   &1 & -1& &  &    \\
   &     &  & & &\ddots     \\
		 \eaa\right)
\quad{\rm and}\quad		 
\A_N=\frac{D}{\dx^2}\left(\baa{cccccc}
  -1&    &   & &  &1   \\
  &\ddots &  & &  &    \\
   & &\ddots& &  &    \\
   1&     &  & & &-1    \\
		 \eaa\right).		
\ee
The exponential of each $\A_j$ can now be evaluated exactly:
\be
\e^{\dt\A_j}=\left(\baa{cccccc}
  1&    &   & &  &   \\
  &\alpha &\beta & &  &    \\
   &\beta &\alpha& &  &    \\
   &     &  & & &1     \\
		 \eaa\right),
\quad\quad		 
\e^{\dt\A_N}=\left(\baa{cccccc}
  \alpha&    &   & &  &\beta   \\
  &1 &  & &  &    \\
   & & &1  &    \\
   \beta&     &  & & &\alpha    \\
		 \eaa\right),		
\ee
where
\be
\alpha=\frac12(1+\gamma),\quad
		 \beta=\frac12(1-\gamma),\quad{\rm and}\quad\gamma=\e^{-2\ld}.
\la{ai}
\ee
Each $\e^{\dt\A_j}$ updates  only $u_j$ and $u_{j+1}$ as
\ba
u_j^\prime&=&\alpha u_j+\beta u_{j+1}\nn\\
u_{j+1}^\prime&=&\beta u_j+\alpha u_{j+1}.
\la{update}
\ea
The eigenvalues of this updating matrix are $\alpha\pm\beta=1,\ \gamma$, with
det $=\gamma$. This means that the updating is dissipative for $r>0$ and
unstable for $r<0$. Since $\alpha$ and $\beta$ are given in terms of $\gamma$, {\it the
resulting algorithm depends only on a single parameter $\gamma$}. 

One can now decompose $\exp(\dt\A)$ to
first order in $\dt$ (apply (\ref{sym1ab}) repeatedly) via either
\be
{T}_{\rm 1A}(\dt)=e^{\dt\A_N}\cdots
e^{\dt\A_2}e^{\dt\A_1},\la{sym1a}
\ee
or
\be
{T}_{\rm 1B}(\dt)=e^{\dt\A_1}\cdots
e^{\dt\A_{N-1}}e^{\dt\A_N}.
\la{sym1b}
\ee
These algorithms update the grid points sequentially, two by two
at a time according to (\ref{update}), but each grid point is updated {\it twice},
in an intertwining manner. 
This is crucial for dealing with the periodic boundary condition.
Let $u^*_j$ denotes 
the first time when $u_j$ is updated and $u^\prime_j$ the second (and final) time it
is updated. One then has
for algorithm 1A:
\ba
u_1^* &=&\alpha u_1+\beta u_2\la{ustart}\\
u_2^* &=&\beta u_1+\alpha u_2\nn\\
u_2^\prime &=&\alpha u_2^*+\beta u_3\nn\\
u_3^* &=&\beta u_2^*+\alpha u_3.\nn\\
&&\cdots\nn\\
u_j^\prime &=&\alpha u_j^*+\beta u_{j+1}\nn\\
u_{j+1}^* &=&\beta u_j^*+\alpha u_{j+1}.\nn\\
&&\cdots\nn\\
u_N^\prime &=&\alpha u_N^*+\beta u_1^*\nn\\
u_1^\prime &=&\beta u_N^*+\alpha u_1^*.
\la{upend}
\ea
Since $\alpha+\beta=1$, summing up both sides from (\ref{ustart}) to (\ref{upend}) gives,
\be
\sum_{j=1}^Nu^\prime_{j}= \sum_{j=1}^Nu_{j}.
\la{np}
\ee  
The algorithm is therefore norm-conserving. 
The same is true of algorithm 1B below.
For $2<j<N$ one has
\ba
u_j^\prime &=&\alpha u_j^*+\beta u_{j+1}\nn\\
&=&\alpha(\beta u^*_{j-1}+\alpha u_j)+\beta u_{j+1}\nn\\
&=&\beta (u_{j-1}^\prime-\beta u_j)+\alpha^2 u_j+\beta u_{j+1}\nn\\
&=&\beta u_{j-1}^\prime+\gamma u_j+\beta u_{j+1}
\la{alg1a}
\ea
and for	$j=2,N$,
\ba
u_2^\prime &=&\beta u_{1}^*+\gamma u_2+\beta u_3\nn\\
u_N^\prime &=&\beta u_{N-1}^\prime+\gamma u_N+\beta u_1^*.
\ea
Finally when the snake bits its tail, one has
\be
u_1^\prime =\frac\beta\alpha u_{N}^\prime+\gamma u_1+\frac{\gamma\beta}\alpha u_2.
\ee
Similarly, 1B is given by
\ba
u_1^* &=&\beta u_N+\alpha u_1\la{dnstar}\\
u_{N}^\prime &=&\beta u_{N-1}+\gamma u_N+\beta u_{1}^*\nn\\
u_j^\prime &=&\beta u_{j-1}+\gamma u_j+\beta u_{j+1}^\prime\la{alg1b}\\
u_{2}^\prime &=&\beta u_{1}^*+\gamma u_2+\beta u_{3}^\prime\nn\\
u_1^\prime &=&\frac{\gamma\beta}\alpha u_{N}+\gamma u_1+\frac{\beta}\alpha u_2^\prime.
\ea

Algorithms 1A and 1B
are essentially given by (\ref{alg1a}) and (\ref{alg1b}) respectively, except for three values of   
$u_1^\prime$, $u_2^\prime$ and $u_N^\prime$. The forms of (\ref{alg1a}) and (\ref{alg1b})
reproduce Saul'yev's schemes (\ref{saula}) and (\ref{saulb}),
but with different coefficients. Saul'yev's coefficient $\gamma\!_S$ is a rational
approximation to the $\gamma=\e^{-2r}$ here. Note that his
$\beta\!_S=r/(1+r)$ is also given by $\beta\!_S=(1-\gamma\!_S)/2$.  
In contrast to Saul'yev's algorithm, which cannont be started for periodic 
boundary condition, algorithm 1A and 1B are truly explicit because they are 
fundamentally given by the sequential updating of (\ref{update}). 
Each algorithm can get started by first updating $u_1$ to $u_1^*$, then updating it 
again at the end to $u^\prime_1$.   

By virtue of (\ref{sym2ab}) one can now
immediately generate a second-order time-marching algorithm via the symmetric product,
\ba
{T}_{2}(\dt)&=&{T}_{\rm 1B}(\frac\dt{2}){T}_{\rm 1A}(\frac\dt{2})\nn\\
&=&e^{\frac12\dt\A_1}e^{\frac12\dt\A_2}
\cdots e^{\frac12\dt\A_N}e^{\frac12\dt\A_N}\cdots
e^{\frac12\dt\A_2}e^{\frac12\dt\A_1}.
\la{symt}
\ea
If the boundary effects of $u_1^\prime$, $u_2^\prime$ and $u_N^\prime$ are ignored (for now)
and 1A and 1B are considered as given by (\ref{alg1a}) and (\ref{alg1b}),  
then the alternative product
${T}_{\rm 1A}(\dt/2){T}_{\rm 1B}(\dt/2)$ yields the same second-order algorithm.
In this case, algorithms 1A and 1B have amplification factors
\ba
g_{\rm 1A}&=&\frac{\gamma+\beta\e^{i\theta}}{1-\beta\e^{-i\theta} },\nn\\
g_{\rm 1B}&=&\frac{\gamma+\beta\e^{-i\theta}}{1-\beta\e^{+i\theta} },
\ea
with opposite phase errors, and the second-order algorithm has 
\ba
g_{2}&=&g_{\rm 1B}\left(\frac{\dt}{2}\right)
g_{\rm 1A}\left(\frac{\dt}{2}\right)\nn\\
&=&\frac{\widetilde\gamma^2
+\widetilde\beta^2+2\widetilde\beta\widetilde\gamma\cos\!\theta}
{1+\widetilde\beta^2-2\widetilde\beta\cos\!\theta},
\la{g2}\\
&=&\frac{1-(4\widetilde\beta\widetilde\gamma/\widetilde\alpha^2)\sin^2\!\theta/2}
{1+(4\widetilde\beta/\widetilde\alpha^2)\sin^2\!\theta/2}=e^{-h_2},
\la{g2a}
\ea
with no phase error and where
\be
\widetilde\alpha=\frac12(1+\widetilde\gamma),\quad\widetilde\beta
=\frac12(1-\widetilde\gamma)\quad{\rm and}\quad\widetilde\gamma=\gamma(r/2).
\la{abi}
\ee
Since both algorithms 1A and 1B have phase errors, only the second-order 
algorithm is qualitatively similar to the exact solution. Eq.(\ref{g2a}) makes it clear that
this algorithm is unconditionally stable since $0\le\widetilde\gamma\le1$. 
Algorithms 1A and 1B are also unconditionally stable since $|g_{\rm 1A,1B}|=\sqrt{g_2}$ with
$\widetilde\gamma\rightarrow\gamma$. 
Note that this also proves the unconditional stability of Saul'yev's algorithms. His coefficient
$\gamma\!_S$ can turn negative, but only approaches -1 as $r\rightarrow\infty$. 
Conventional explicit methods, like that of the Euler algorithm, are limited in stability
because they have power-series amplification factors. By contrast, symplectic
finite-difference methods are unconditionally stable because they produce Saul'yev-type 
schemes with {\it rational-function} amplification factors. 
This type of stability is usually associated only with {\it implicit} methods.

\begin{figure}[hbt]
\includegraphics[width=0.70\linewidth]{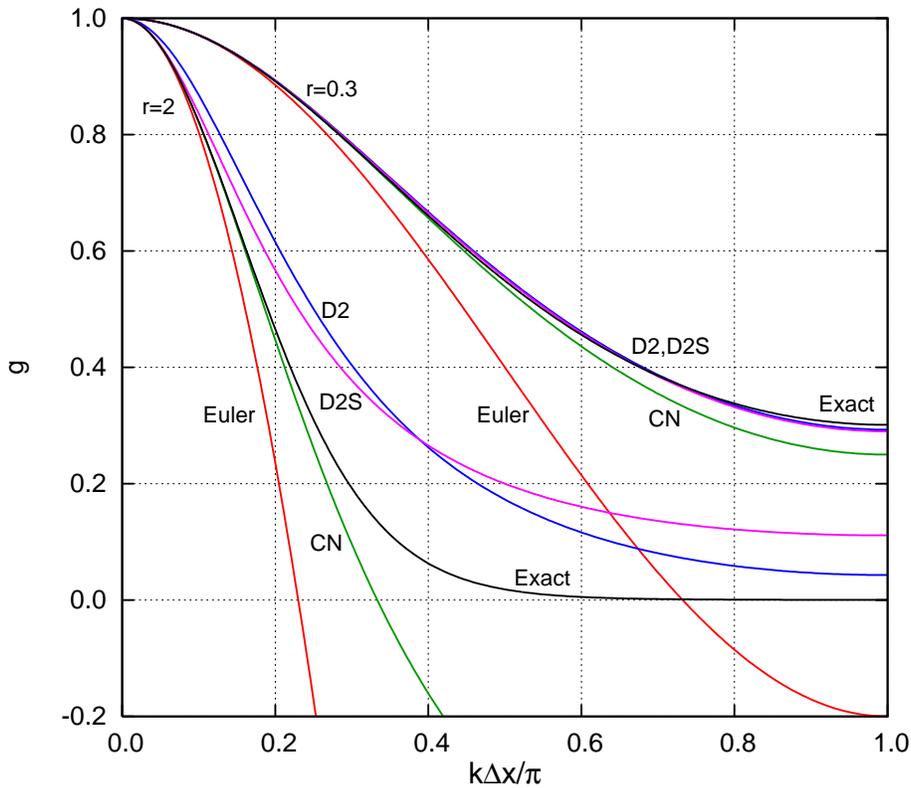}
\caption[]{\label{stable} Comparing the amplification factor of various diffusion
algorithms at $r=0.3$ and $r=2$. Euler is the first-order explicit algorithm. 
CN is the second-order implicit Crank-Nicolson algorithm. D2 is the second-order 
symplectic algorithm using the original coefficients (\ref{abi}) and D2S is the 
same alogrithm but uses Saul'pev's coefficient (\ref{saulcoef}).
``Exact" is $g_{ex}$ of (\ref{gex}). 
}
\end{figure}

Since the order of matrices defining $T_2(\dt)$ in (\ref{symt}) is left-right symmetric, one has
the same situation as in the symplectic integrators case of (\ref{two1}), implying 
that $g_2(-r)g_2(r)=1$. 
This means that $g_2(-r)> 1$ for all $k\ne 0$ modes and 
the algorithm is unconditionally unstable for negative time steps. 
Moreover, this also means that $h_2(-r)=-h_2(r)$ and $h_2(r)$ is an odd function of $r$. 
This is a simpler, functional version of (\ref{oddstr}). 
Expanding $h_2$ of (\ref{g2a}) in powers of $\theta$ gives,
\be
h_2=\frac{2(1-\widetilde\gamma)}{1+\widetilde\gamma} \theta^2
-\frac{(13-35\widetilde\gamma+35\widetilde\gamma^2-13\widetilde\gamma^3) }{6(1+\widetilde\gamma)^3}\theta^4+\cdots .
\la{g22}
\ee
Each coefficient must be an odd function of $r$. This is satisfied only if
\be
\widetilde\gamma(-r)\widetilde\gamma(r)=1.
\la{gcond}
\ee
Thus every function $\widetilde\gamma(r)$ satisfying (\ref{gcond}) with $|\widetilde\gamma(r)|\le 1$
defines an unconditionally stable algorithm for solving the diffusion equation. 
For the original choice of $\widetilde\gamma(r)=\e^{-r}$, one finds
\be
h_2=(r-\frac{r^3}{12}+\frac{r^5}{120}\cdots) \theta^2
-(\frac{r}{12}+\frac{35 r^3}{144}+\cdots)\theta^4
+(\frac{r}{360}+\frac{539r^3}{4320}+\cdots)\theta^6 +\cdots .
\la{h2org}
\ee
Comparing this to the expansion of the 
exact amplification exponent,
\be
h_{ex}=r \theta^2-\frac{r}{12}\theta^4+\frac{r}{360}\theta^6 +\cdots,
\la{gex2}
\ee
one sees that the original choice does not 
reproduce leading term exactly except when $r<<1$. 
To improve this, let's take $\widetilde\gamma(r)$ to be an arbitrary function of $r$
but with $\widetilde\alpha$ and $\widetilde\beta$ still defined by (\ref{abi}). 
The first term in $h_{ex}$ can now be matched exactly by requiring
\be
\frac{1-\widetilde\gamma(r)}{1+\widetilde\gamma(r)}=\frac{r}2
\quad\rightarrow\quad
\widetilde\gamma(r)=\frac{1-r/2}{1+r/2},
\la{saulcoef}
\ee
which is precisely Saul'yev's original coefficient. With this choice for $\widetilde\gamma(r)$, (\ref{g2a}) reads
\be
g_2=\frac{1-2r(1-r/2)\sin^2(\theta/2)}
{1+2r(1+r/2)\sin^2(\theta/2)},
\la{g23}
\ee
with exponent
\be
h_2=r \theta^2-(\frac{r}{12}+\frac{r^3}{4})\theta^4+(\frac{r}{360}+\frac{r^3}{8}+\frac{r^5}{16})\theta^6 +\cdots 
\la{h2sp}
\ee
which is now correct to third-order in $\theta$. 
Comparing this and (\ref{h2org}) to $h_{ex}$, one sees that all the error terms of $h_2$
are {\it odd} powers of $r$ higher than the first. As a matter of fact, by resumming terms proportional to $r$,
we can make this error structure in exact conformity with (\ref{oddstr}),
\be
h_2=rh_0+r^3E_3(\theta)+r^5E_5(\theta)+\cdots,
\la{erst}
\ee
where $h_0=4\sin^2(\theta/2)$.  
This is the same error structure exploited by symplectic integrators to produce higher order algorithms.
 
Saul'yev's algorithm is close in reproducing the 
amplification factor of the {\it implicit} Crank-Nicolson (CN) scheme 
(which is without the $\pm r/2$ terms in (\ref{g23})).
The CN scheme has the advantage that its exponent is
\be
h_{CN}=r \theta^2-\frac{r}{12}\theta^4+(\frac{r}{360}+\frac{r^3}{12})\theta^6 +\cdots 
\ee
which is correct to fifth-order in $\theta$. 
In Fig.\ref{stable}, we compare the amplification factor of various algorithms at two values of $r$.
D2 and D2S are second-order symplectic diffusion algorithms described above using the original coefficient (\ref{abi})
and Saul'pev's coefficient (\ref{saulcoef}), respectively. 
For $r=0.3$, both D2 and D2S track $g_{ex}$ closely over the entire range $k$ values. Both Euler 
and CN tend to over-damp higher Fourier modes. At $r=2.0$, while both Euler and CN turn 
negative at large $k$, D2 and D2S remain positive, like that of $g_{ex}$. 
At large $r$, D2S is clearly better than D2 at small $k$.  

\begin{figure}[hbt]
\includegraphics[width=0.70\linewidth]{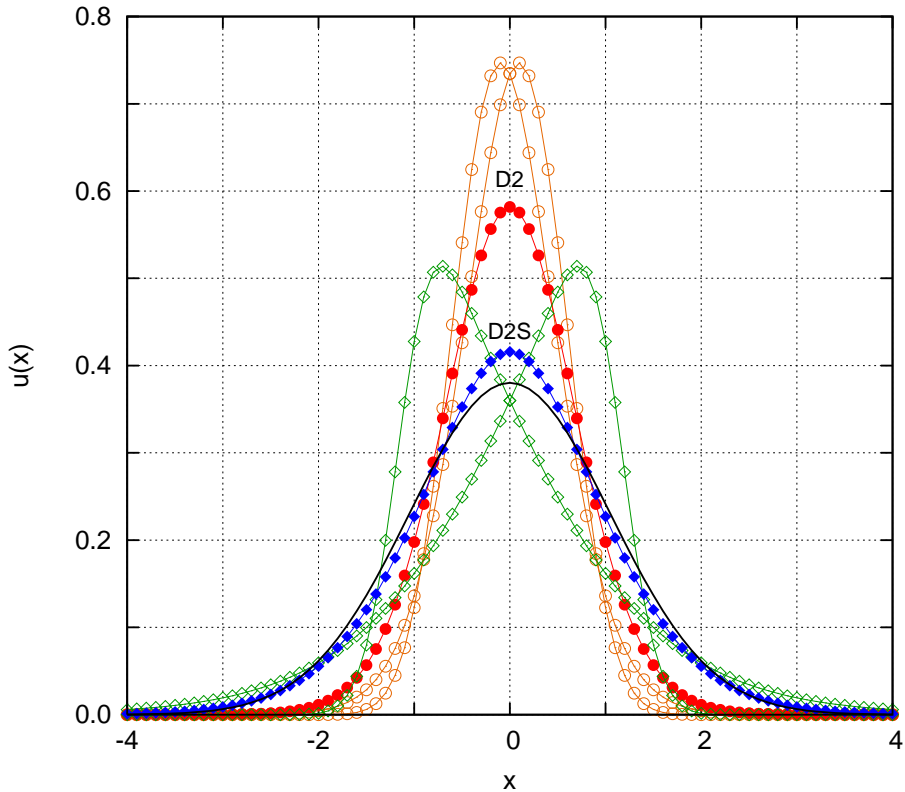}
\caption[]{\label{diffuse} The diffusion of a Gaussian profile after $t=1$ on a grid of 120 points
spanning the interval [-6,6] with $D=1/2$ and $\dt=0.1$, corresponding to $r=5$. This large value is
choosen to exaggerate various errors. The open and filled circles are algorithms 1A, 1B and 2 
using the standard coefficients. The open and filled diamonds are the same three algorithms 
using Saul'pev's coefficients. The solid black line is the exact solution in the continuum limit. 
}
\end{figure}

To generalize Saul'yev schemes to periodic boundary condition, one simply replaces in the above
algorithms, $\gamma\rightarrow\gamma\!_S$. In Fig.\ref{diffuse} we show the working of algorthms
1A, 1B and 2 using both sets of cofficients. With the original coefficients, the 
phase errors are much smaller, but the under-damp error is much larger. 
For Saul'pev's coefficient, the phase errors are much greater, but the under-damp error 
is smaller. Since the phase error is automatically eliminated by going to second-order,
D2S has an advantage over D2.   

If one were to construct a 
2n$^{th}$-order algorithm out of a product of second-order algorithms
\be
T_{2n}(\dt)=\prod_{i=1}^{N} T_2(a_i\dt),
\la{t2n2}
\ee
then the corresponding $h_{2n}$ is given by
\ba
h_{2n}(r)&=&\sum_{i=1}^Nh_2(a_ir)\nn\\
&=&r h_0\sum_{i=1}^Na_i+r^3E_3(\theta)\sum_{i=1}^Na^3_i+r^5E_5(\theta)\sum_{i=1}^Na^5_i+\cdots,
\ea
and the order conditions (\ref{sum1})-(\ref{sum5}) are easily understood.
Unfortunately, for the diffusion algorithm, $T_2(a_i\dt)$ is unstable for
any negative $a_i$ and no negative coefficient $a_i$ can be allowed. In this case,
order-conditions such as (\ref{sum3}) and (\ref{sum5}) cannot be satisfied and no higher-order 
composition algorithms of the form (\ref{t2n2}) is possible.
(However, these algorithms will be useful for solving the advection equation in the next section.)

To overcome this impass, one must go beyond the single product approximation of (\ref{t2n2}), and
consider a multi-product expansion of the form
\be
T_{2n}(\dt)=\sum_k c_k\prod_{i=1}^{N_k} T_2(a_{ki}\dt).
\la{mt2n}
\ee
If $a_{ki}$ were to remain positive, then Shang\cite{sheng89} has shown that
any single product in (\ref{mt2n}) can at most be second order and that 
some ${c_k}$ coefficients must be negative.  
More recently this author realize that\cite{chin10,chin11}, due to the error structure (\ref{erst}),
the first-order term $rh_0$ is automatically preserved by monomial products of the form
$T_2^k(\dt/k)$ with	exponent
\be
h_{2}^{(k)}=rh_0+k^{-2}r^3E_3(\theta)+k^{-4}r^5E_5(\theta)+\cdots.
\ee
The arbitrariness in $N_k$ and $a_{ki}$ can be eliminated by taking $N_k=k$ and $a_{ki}=1/k$.
This then produces a much simpler Multi-Product Expansion\cite{chin10} (MPE)
\be
T_{2n}(\dt)=\sum_k c_k T_2^k(\dt/k)
\ee
for any sequence of $n$ whole numbers $\{k\}$ with analytically known coefficients $c_k$.
For the harmonic sequence of $k=1,2,3, \cdots$,
the first few higher order algorithms are:
\be
{T}_4(\ep)=-\frac13{T}_2(\ep)
+\frac43{T}_2^2\left(\frac\ep{2}\right)
\la{four}
\ee
\be
{T}_6(\ep)=\frac1{24} {T}_2(\ep)
-\frac{16}{15}{T}_2^2\left(\frac\ep{2}\right)
+\frac{81}{40}{T}_2^3\left(\frac\ep{3}\right)
\la{six}
\ee
\be
{T}_8(\ep)=-\frac1{360} {T}_2(\ep)
+\frac{16}{45}{T}_2^2\left(\frac\ep{2}\right)
-\frac{729}{280}{T}_2^3\left(\frac\ep{3}\right)
+\frac{1024}{315}{T}_2^4\left(\frac\ep{4}\right).
\la{eight}						   
\ee
For the diffusion equation, the use of the fourth-order extrapolation (\ref{four}) 
has been previously suggested by Schatzman\cite{schatz02}. However, these high 
order methods are uselss for conventional explicit schemes, since they are 
unstable at large time steps. It is only with the use of unconditionally stable 
algorithms here that the power of these high order schemes can be unleashed. 
These MPE algorithms do not preserve the positivity of the 
initial profile. This is in keeping with the general observation that there 
can't be any finite-difference scheme for solving the diffusion equation that 
preserves positivity beyond the second-order\cite{hv}.
More recently, Zillich, Mayrhofer and Chin\cite{zmc10} have shown that
Path-Integral Monte Carlo simulations, where positivity is of the utmost importance, 
can be successfully carried out using these expansions, demonstrating that the violation 
of positivity is small and controllable. 
These MPE algorithms are also not symplectic. However, as argued by Blanes, Casas and Ros\cite{bcr99}, 
they are symplectic to order $2n+3$. Thus at sufficiently high orders, they are 
indistinguishable from truly symplectic algorithms up to machine precision.
In the context of classical dynamics, these MPE algorithms have been tested up 
to the 100$^{th}$ order in Ref.\onlinecite{chin11}.

\begin{figure}[hbt]
\includegraphics[width=0.70\linewidth]{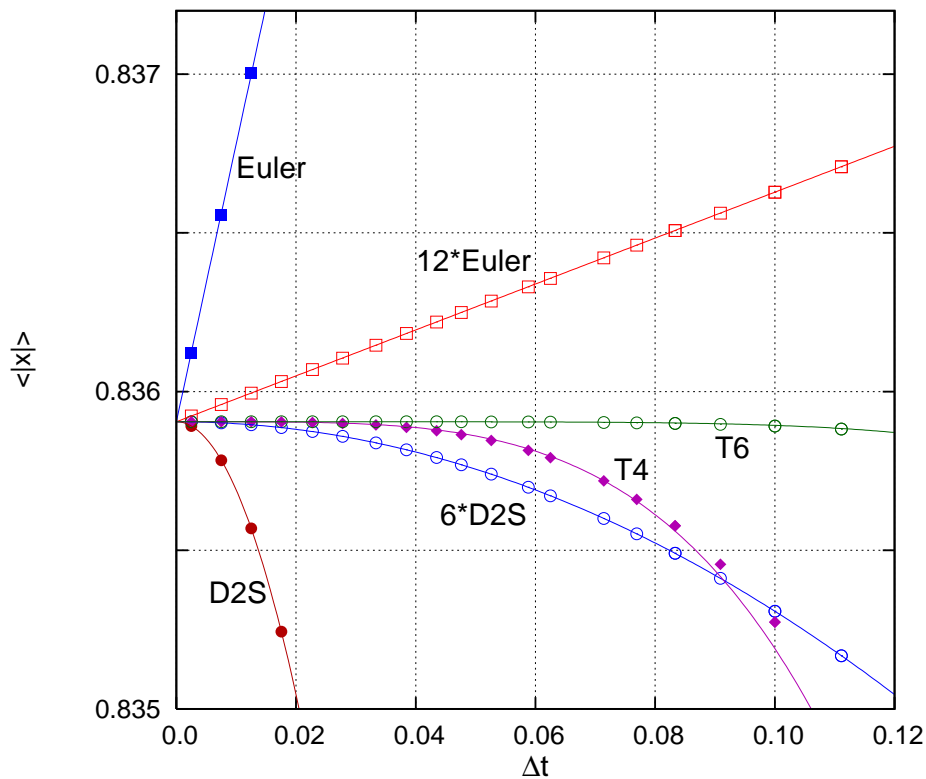}
\caption[]{\label{cong} The convergence of $\langle |x|\rangle$ after a Gaussian profile 
has been diffused for $t=1$ on a grid of 120 points
spanning the interval [-6,6] with $D=1/2$. The range of $\dt$, 
from $\dt=1/400$ to $\dt=1/9$, corresponds to a range of $r=0.125$ to $r=5.555$.
Lines are fitted power laws $\dt^n$ verifying the order of
the algorithm. ``12Euler" labels results of running the first order Euler algorithm 12 times at
time step $\dt/12$. ``6D2S" are results from running the second-order algorithm D2S six times at
time step $\dt/6$. $T_4$ and $T_6$ are fourth and sixth-order algoirthms which 
require 3 and 6 runs of D2S respectively.  
}
\end{figure}

In Fig.\ref{cong}, we compare and verify the order of
convergence of $T_4$ and $T_6$ with D2S as $T_2$ by
computing the expectation value
\be
\langle|x|\rangle=\frac{\sum_{j=1}^N|j\dx|u_j}{\sum_{j=1}^N u_j}
\ee
after evolving $u_j$ for $t=1$ as a function of $\dt$.
The absolute value is used because the Euler algorithm would exactly preserve $\langle x\rangle$
even when it is unstable. For computing $\langle|x|\rangle$,
the Euler algorithm is extremely linear within its tiny range of
stability. Since it tends to over-damp, its evolving profile is flatter and $\langle|x|\rangle$ 
converges from above. Symplectic diffusion algorithms under-damp, and their results 
for $\langle|x|\rangle$ converge from below. All results can be well fitted with power 
laws of the form $a+b\dt^n$ with $n=1,2,4$ and 6, verifying the order of the algorithms.
To show that these higher order algorithms are more efficient than running low order algorithms 
at reduced step sizes, we also plotted results of running the Euler algorithm
12 times at time step $\dt/12$ and algorithm D2S 
six times at $\dt/6$.

\section {Symplectic advection algorithm}

For the advection equation
\be
\frac{\partial u}{\partial t}=-v\frac{\partial u}{\partial x},
\la{ad}
\ee
its usual semi-discrete form is
\be
\frac{\partial u_j}{\partial t}=-\frac{v}{2\dx}(u_{j+1}-u_{j-1}),
\ee
with discretization matrix
\be
\B=\frac{v}{2\dx}
\left(\baa{ccccc}
0 &  -1  &   & &  1   \\
 1 &  0 & -1 & &      \\
   & & \ddots &  &     \\
   &     &   1& 0& -1    \\
 -1  &     &   & 1& 0    \\
\eaa\right),
\ee
and solution
\be
\bu(t+\dt)=e^{\dt\B}\bu(t).
\la{beq}
\ee
The exact amplification factor ($\theta=k\dx$)
\be
g_{ex}=\e^{-i\eta\sin\theta}, \qquad {\rm with}\qquad \eta=\frac{v\dt}{\dx},
\la{gad}
\ee
is unitary and causes a phase-shift of each Fourier component. In the
limit of $\dx\rightarrow0$, the phase-shift becomes uniform for all the Fourier
components $\e^{ikx}\rightarrow\e^{ik(x-v\dt)}$, resulting in a uniform shift of
the entire function $u(x)\rightarrow u(x-v\dt)$, which is the exact solution to
(\ref{ad}). Any Taylor expansion of (\ref{beq}) will produce algorithms with a non-unitary $g$,
resulting in unwanted dissipations or instability. The situation here is much more delicate 
than in the diffusion case. 

The natural decomposition is similarly,
\be
\B=\sum_{i=1}^N\B_j,
\ee
where 
\be
\B_j=\frac{v}{2\dx}\left(\baa{cccccc}
  \ddots&    &   & &  &   \\
  & 0 & -1 & &  &    \\
   &1 & 0& &  &    \\
   &     &  & & &\ddots     \\
		 \eaa\right)
\quad{\rm with}\quad		 
\B_N=\frac{v}{2\dx}\left(\baa{cccccc}
  0&    &   & &  &1   \\
  &\ddots &  & &  &    \\
   & &\ddots& &  &    \\
   -1&     &  & & &0    \\
		 \eaa\right).		
\ee
It follows that
\be
\e^{\dt\B_j}=\left(\baa{cccccc}
  1&    &   & &  &   \\
  &c &-s & &  &    \\
   &s &c& &  &    \\
   &     &  & & &1     \\
		 \eaa\right),
\quad\quad		 
\e^{\dt\B_N}=\left(\baa{cccccc}
  c&    &   & &  &s   \\
  &1 &  & &  &    \\
   & & &1  &    \\
   -s&     &  & & &c    \\
		 \eaa\right),		
\ee
where now
\be
c=\cos(\eta/2)\quad{\rm and}\quad s=\sin(\eta/2).
\la{adv}
\ee
Each $\e^{\dt\B_j}$ only updates $u_j$ and $u_{j+1}$ as
\ba
u_j^\prime&=& cu_j-su_{j+1}\nn\\
u_{j+1}^\prime&=&su_j+cu_{j+1}.
\la{adupm}
\ea 
As in the diffusion case, the above updating can be recasted into the following 
forms for algorithms 1A and 1B, with 1A given by
\ba
u_1^* &=&c u_1-s u_2\la{adustart}\nn\\
u_2^\prime &=& su_{1}^*+ u_2-s u_3\nn\\
u_j^\prime 
&=&su_{j-1}^\prime+u_j-s u_{j+1} \quad (2<j<N)
\la{ad1a}\\
u_N^\prime &=& su_{N-1}^\prime+ u_N-su_1^*\nn\\
cu_1^\prime &=& s u_{N}^\prime+cu_1-s u_2.\nn
\nn
\ea
and 1B given by
\ba
u_1^* &=&s u_N+c u_1\nn\\
u_{N}^\prime &=&s u_{N-1}+ u_N-s u_{1}^*\nn\\
u_j^\prime
&=&su_{j-1}+u_j-s u_{j+1}^\prime\quad (2<j<N)
\la{ad1b}\\
u_{2}^\prime &=&s u_{1}^*+ u_2-s u_{3}^\prime\nn\\
cu_1^\prime &=& s u_{N}+cu_1-s u_2^\prime.\nn
\ea
In contrast to the diffusion case, these algorithms are {\it not} exactly norm-preserving
for periodic boundary condition. By adding up both sides of the above algorithms, one finds
that what is preserved by 1A is not the usual norm
$
N=\sum_{j=1}^Nu_{j}
$,
but a modified norm given by 
\be
\tilde N_{\rm 1A}=N+(\frac{c}{1-s}-1)u_1
\la{n1a}
\ee
Similarly, what is preserved by 1B is 
\be
\tilde N_{\rm 1B}=N+(\frac{c}{1+s}-1)u_1.
\la{n1b}
\ee
If initially $u_1=0$, then $\tilde N_{\rm 1A}=\tilde N_{\rm 1B}=N_0$, where $N_0$ is the initial norm. 
As the system evolves, each algorithm's actual norm will evolve as
\ba
N_{\rm 1A}&=&N_0-(\frac{c}{1-s}-1)u_1,\nn\\
N_{\rm 1B}&=&N_0-(\frac{c}{1+s}-1)u_1.
\ea   
The error is due to a single point $u_1$, where it is the only point not updated twice immediately.
As the wave form travels around the periodic box, $u_1$ will trace out the shape of the wave and imprint
that as the error of the norm in time. For a sharp pulse, the norm error will return to zero after 
the pulse peak has passed through $u_1$. Thus norm-preservation will be periodic. 
For the advection equation, this is a small effect, and is secondary to the phase and oscillation 
error mentioned below. However, this error will be important in the next section. 

If the boundary values $u_1^\prime$, $u_2^\prime$ and $u_N^\prime$ are ignored for now, then 
again the resulting second-order algorithm is unique, independent of the order of
applying 1A or 1B. The amplification factors are all {\it unitary}:
\ba
g_{\rm 1A}&=&\frac{1-s\e^{i\theta}}{1-s\e^{-i\theta} }=\exp(-i\phi_{\rm 1A}),\la{gad1a}\\
g_{\rm 1B}&=&\frac{1+s\e^{-i\theta}}{1+s\e^{i\theta} }=\exp(-i\phi_{\rm 1B},)\la{gad1b}\\
g_{2}&=&g_{\rm 1B}(\dt/2)g_{\rm 1A}(\dt/2)\nn\\
&=&\frac{1-i2(\tilde s/\tilde c^2)\sin\theta}{1+i2(\tilde s/\tilde c^2)\sin\theta }=\exp(-i\phi_2),
\la{gad2}
\ea
with phase angles
\ba
\phi_{\rm 1A}
&=&2\tan^{-1}\left(\frac{s\sin\theta}{1-s\cos\theta}\right),\nn\\
\phi_{\rm 1B}
&=&2\tan^{-1}\left(\frac{s\sin\theta}{1+s\cos\theta}\right),\nn\\
\phi_{2}&=&\phi_{\rm 1A}(\dt/2)+\phi_{\rm 1B}(\dt/2),\nn\\
&=&2\tan^{-1}\left( \frac{2\tilde s}{1-\tilde s^2}\sin\theta\right),
\ea	   
where here
\be
\tilde s=\sin(\eta/4)\quad{\rm and}\quad \tilde c=\cos(\eta/4).
\la{adorig}
\ee
Since $g_{\rm 1A}$ and $g_{\rm 1B}$ are not complex conjugate of each other,
their phase errors do not exactly cancel. Their residual difference is the error of
the second-order algorithm. 

Algorithms (\ref{ad1a}) and (\ref{ad1b}) are the corresponding Saul'yev's schemes
for solving the advection equation. The coefficient here
is $s=\sin(\eta/2)$ rather than Saul'yev's coefficient of $s=\eta/2$. This explains
why it makes no sense to apply Saul'yev's schemes at $s>1$, since they can no longer be
derived from the fundamental updating matrix (\ref{adupm}) with a real $c=\sqrt{1-s^2}$.
At $s>1$, Saul'yev's schemes are in fact unstable, suffering from spatial amplification\cite{cy07},
despite the unimodulus appearance of (\ref{gad1a}) and (\ref{gad1b}). This is easy to see in the
case of algorithm 1A. If initially $u_j=0$ for $j\ge J$, but $u_{J-1}\ne 0$,  
then according to (\ref{ad1a}), $u^\prime_{J+n}=s^{n+1}u^\prime_{J-1}$ increases without bound
as a function of $n$. Even the case of $s=1$ is pathological. For Saul'yev's coefficient $s=\eta/2=1$, 
one has
\be
g_{\rm 1A}=-\e^{ik\dx}=-\e^{ik(v/2)\dt}\quad{\rm and}\quad g_{\rm 1B}=\e^{-ik\dx}=\e^{-ik(v/2)\dt}.
\ee
Under algorithm 1A, Fourier mode $\e^{ikx}$ will flip its sign and propagate with
velocity $-v/2$. Under 1B, it will propagate with velocity $v/2$. The resulting second
order algorithm then leaves the Fourier mode {\it stationary} with only a sign flip. 
This is completely contrary to the behavior of the exact solution and is a source of great 
error for Saul'yev's schemes. As we will show below, alternative choices for $s$ will
eliminate such unphysical behaviors.       

While the derived choice of $s=\sin(\eta/2)$ is unconditionally stable for
all $\eta$, the resulting algorithms 1A and 1B have {\it huge} phase errors, and are no better
than Saul'yev's choice of $s=\eta/2$. This is because in comparison with the exact phase angle,
\be
\phi_{ex}=\eta\sin\theta=\eta\theta-\frac\eta{6}\theta^3+\cdots
\ee 
algorithms 1A and 1B have expansions
\ba
\phi_{\rm 1A}
&=&\frac{2s}{1-s}\theta-\frac{s(1+s)}{3(1-s)^3}\theta^3+\cdots,\nn\\
\phi_{\rm 1B}
&=&\frac{2s}{1+s}\theta-\frac{s(1-s)}{3(1+s)^3}\theta^3+\cdots,\nn\\
\ea	
and neither $s=\eta/2$ nor $s=\sin(\eta/2)$ can result in a first-order coefficient of $\theta$ 
matching that of $\phi_{ex}$ exactly. The choices of $s$ that can do this are, for 1A,
\be
\frac{2s}{1-s}=\eta \quad\rightarrow\quad s=\frac{\eta}{2+\eta},
\la{ch1a}
\ee
and for 1B,
\be
\frac{2s}{1+s}=\eta \quad\rightarrow\quad s=\frac{\eta}{2-\eta}.
\la{ch1b}
\ee
This then reproduces the Roberts and Weiss\cite{rw66,cy07} forms of the Saul'yev-type algorithm
and will be denoted as RW1A and RW1B. 
For $\eta>0$,
only RW1A is unconditionally stable and RW1B is limited by spatial amplification to $\eta<1$.
The pathological behavior of 1A at $s=1$ can no longer occur at any finite $\eta$.  

For the above choices of $s$, the corresponding phase angles are
\ba
\phi_{\rm 1A}
&=&\eta\theta-(\frac\eta{6}+\frac{\eta^2}{4}+\frac{\eta^3}{12} )\theta^3
+\cdots, \nn\\
\phi_{\rm 1B}
&=&\eta\theta-(\frac\eta{6}-\frac{\eta^2}{4}+\frac{\eta^3}{12} )\theta^3
+\cdots,
\ea
and the modified norms are
\ba
\tilde N_{\rm 1A}&=&N+(\sqrt{1+\eta}-1)u_1\la{nrw1a},\\
\tilde N_{\rm 1B}&=&N+(\sqrt{1-\eta}-1)u_1.
\la{nrw1b}
\ea

The second order algorithm from concatenating RW1A and RW1B is 
\ba
\phi_2&=&\phi_{\rm 1A}(\eta/2)+ \phi_{\rm 1B}(\eta/2)\nn\\
     &=&\eta\theta-(\frac\eta{6}+\frac{\eta^3}{48})\theta^3+(\frac\eta{120}+\frac{5\eta^3}{192}
+\frac{\eta^5}{1280})\theta^5+\cdots
\ea
This second-order advection algorithm will be denoted as RW2.
Because RW1B is limited by spatial amplification to $\eta<1$, RW2 is limited in
stability to $\eta<2$. 	

To generate a stable second-order algorithm for all $\eta$, one can concatenate 1A and 1B
with the same $\tilde s$. To match $g_{ex}$ to first order in $\theta$ then requires  
\be
\frac{2\tilde s}{1-\tilde s}+\frac{2\tilde s}{1-\tilde s}=\eta
\quad\rightarrow\quad
\frac{2\tilde s}{1-\tilde s^2}=\frac\eta{2}\quad\rightarrow\quad 
\tilde s=\frac2{\eta}\left(\sqrt{1+\frac{\eta^2}4 } -1\right).
\la{stilde}
\ee
The resulting amplification factor is, according to (\ref{gad2}),
\be																		  
g_2=\frac{1-i(\eta/2)\sin\theta}{1+i(\eta/2)\sin\theta },
\la{g2cn}
\ee   
which is precisely the {\it implicit} Crank-Nicolson amplification factor.
Since by (\ref{stilde}), $\tilde s\le 1 $, and $\tilde c=\sqrt{1-\tilde s^2}$ is well-defined for all $\eta$,
the algorithm is unconditionally stable and can be applied to periodic boundary problems via the the 
fundamental updating (\ref{adupm}).
Corresponding to (\ref{g2cn}), the phase-angle has the characteristic expansion,
\ba
\phi_2&=&\eta\theta-(\frac\eta{6}+\frac{\eta^3}{12})\theta^3+(\frac\eta{120}+\frac{\eta^3}{24}
+\frac{3\eta^5}{240})\theta^5+\cdots\la{phi2cn}\\
&=&\eta\sin\theta+\eta^3F_3(\theta)+\eta^5F_5(\theta)+\cdots
\ea
where now the time parameter is $\eta$ and the original ``Hamiltonian" is $h_0=\sin\theta$.
We shall designate this second-order algorithm, with $\tilde s$ given by (\ref{stilde}), as A2C.
The second-order algorithm corresponding to 
Saul'pev's choice of $\tilde s=\eta/4$ will be denoted as A2S, and the initially derived
result of $\tilde s=\sin(\eta/4)$ as A2.

\begin{figure}[hbt]
\includegraphics[width=0.49\linewidth]{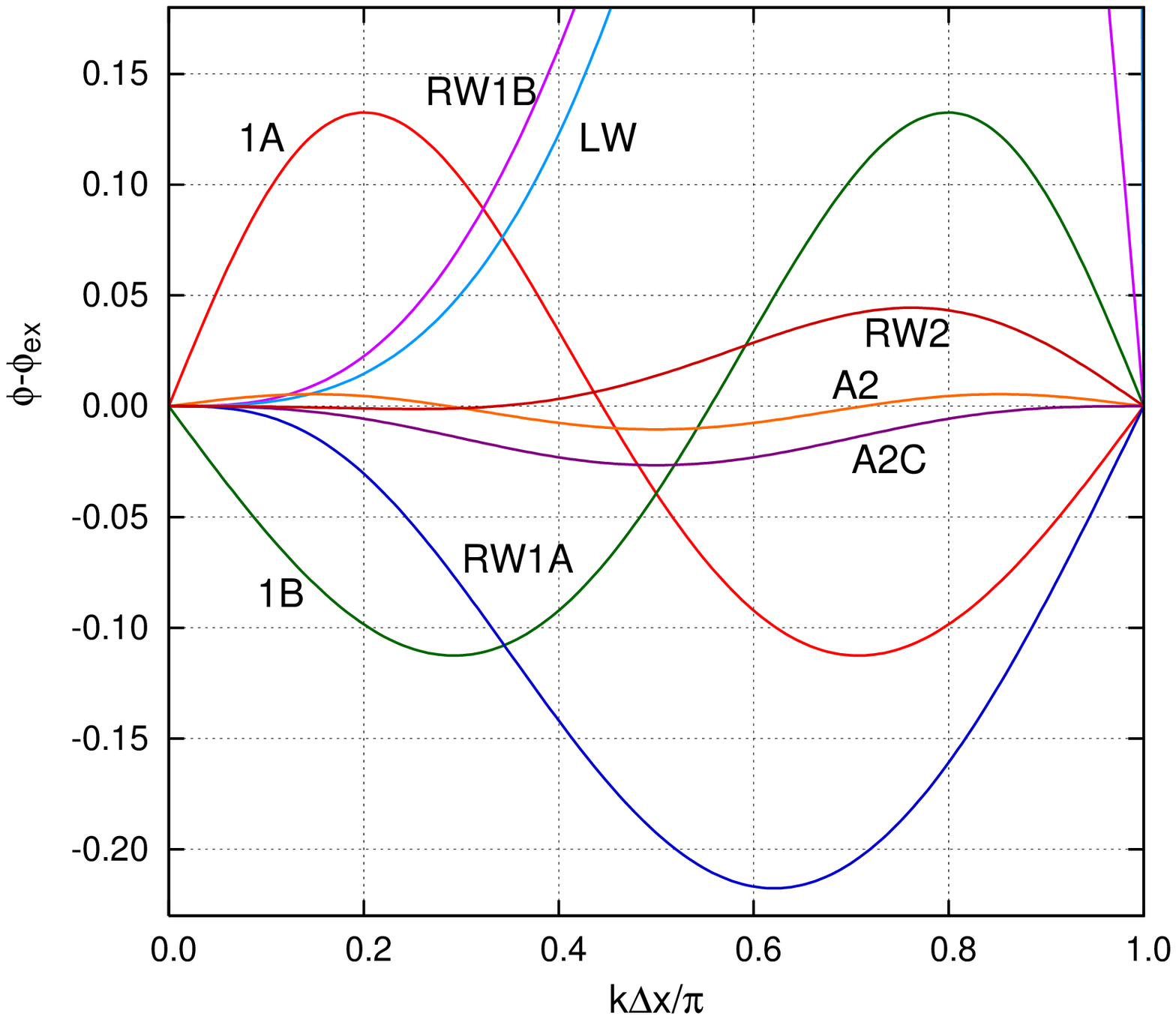}
\includegraphics[width=0.49\linewidth]{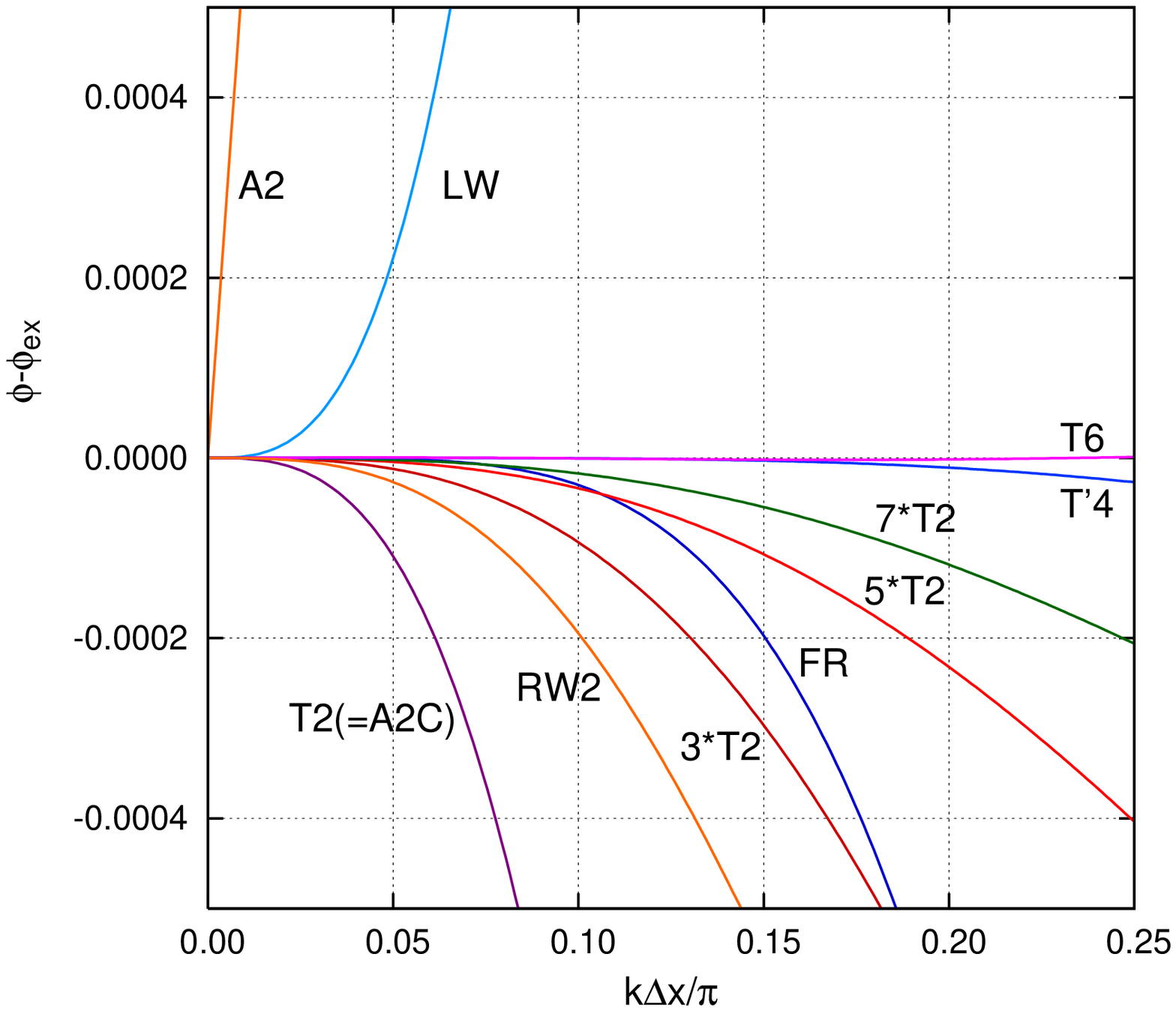}
\caption[]{\label{stbad} 
The phase error of various advection algorithms at $\eta=0.7$. {\bf Left}: The unmodified symplectic
algorithms are denoted as 1A, 1B, and A2. The Robert-Weiss versions are denoted as
RW1A, RW1B and RW2. A2C is the second-order algorithm with the Crank-Nicolson amplification
factor. LW is the Lax-Wendroff scheme included for comparison. {\bf Right:} The phase errors
of fourth and sixth-order algorithms composed out of  
three, five and seven second-order algorithms $T_2$. Their phase errors are compared to that of running 
the second-order algorithm three, five and seven times at reduced time steps. The $T_2$ used here is A2C.
Previous second-order algorithms are also included for a close-up comparison. 
}
\end{figure}

On the left of Fig.\ref{stbad}, the phase error of these symplectic algorithms are exaggerated and compared to
the explicit but dissipative Lax-Wendroff (LW) scheme at a large value of $\eta=0.7$.
The original 1A and 1B algorithms have huge phase errors but are mostly cancelled in the second-order
algorithm A2. Even so, algorithm A2's error curve has a finite slope at $\theta=0$, as shown on the right of
 Fig.\ref{stbad}. By construction, schemes RW1A, RW1B, RW2, A2C and LW all have
zero error slopes at $\theta=0$. This is a crucial advantage of A2C over A2. Also,
A2C is stable for all $\eta$, while RW2 is limited by spatial amplification to $\eta<2$.

The importance of having a zero error slope in the phase angle can be better appreciated from
the following considerations. Let
\be
\chi(t)=\frac{\int x u(x,t)dx}{\int u(x,t)dx}.
\ee 
Since the solution to the advection equation is $u(x,t)=u_0(x-vt)$, we have the
exact result
\be
\chi(\dt)=\frac{\int x u_0(x-v\dt)dx}{\int u_0(x-v\dt)dx}
=\frac{\int (y+v\dt) u_0(y)dy}{\int u_0(y)dy}=\chi(0)+v\dt.
\la{advec}
\ee 
Multiplying algorithm 1A (\ref{ad1a}) by $j$ and sum over $j$ 
yields 
\ba
\sum_{j=1}^Nj u_j^\prime-s\sum_{j=1}^Nju_{j-1}^\prime
&=&\sum_{j=1}^Nju_j-s\sum_{j=1}^Nju_{j+1}\nn\\
(1-s)\sum_{j=1}^Nj u_j^\prime-s\sum_{j=1}^Nu_{j-1}^\prime
&=&(1-s)\sum_{j=1}^Nju_j+s\sum_{j=1}^Nu_{j+1}.
\ea
For a localized pulse far from the boundary, the norm can be consider conserved, 
\be
\sum_{j=1}^Nu_{j-1}^\prime=\sum_{j=1}^Nu_{j+1}.
\ee
One then has the discrete version of (\ref{advec}) 
\be
\frac{\sum_{j=1}^N(j\dx) u_j^\prime}{\sum_{j=1}^Nu_{j}^\prime}
=\frac{\sum_{j=1}^N(j\dx)u_j}{\sum_{j=1}^Nu_{j}}+\frac{2s}{1-s}\dx
\ee
which will reproduce the displacement exactly if
\be
\frac{2s}{1-s}=\frac{v\dt}\dx=\eta,
\ee
which is the condition (\ref{ch1a}) for a zero error-slope.
Similarly for 1B satisfying (\ref{ch1b}). Far from the boundary,
symplectic advection algorithms with 
a zero-error slope in the phase angle would exactly preserve the first two moments of 
$\langle x^n \rangle$. 

\begin{figure}[hbt]
\includegraphics[width=0.7\linewidth]{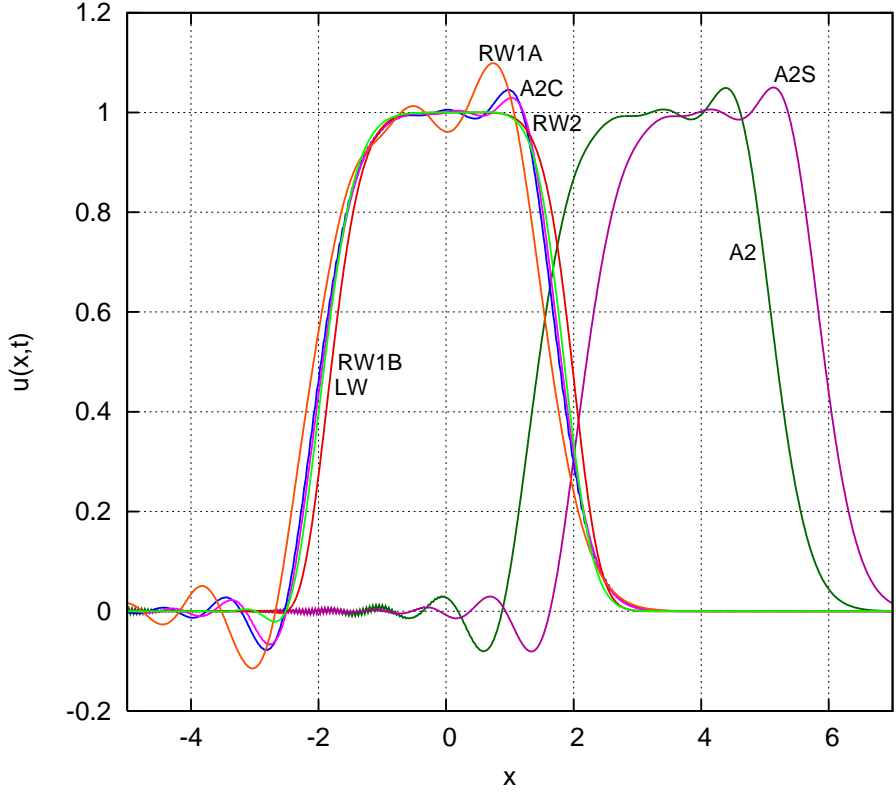}
\caption[]{\label{alg7} 
The propagation of initial profile (\ref{intp})
fives times around a periodic box of [-10,10] with $\dx=0.025$, $\dt=0.02$, $v=1$, and $\eta=0.8$,
corresponding to 5000 iterations of each algorithm. If there were no phase error, the profile would
remain centered on $x=0$. Algorithms A2 and A2S have large and positive, phase errors. 
The oscillation errors in these second-order symplectic algorithms are 
predominately due to the imbedded 1A algorithm. Algorithm RW1B has a much smaller oscillation error 
and is comparable to the oscillation error in the dissipative Lax-Wendroff scheme 
(bright green line).  
}
\end{figure}

In Fig.\ref{alg7} we show the working of these algorithms in propagating an initial profile
\be
u(x,0)=\exp\left[-\left( \frac{x}2 \right)^6\right].
\la{intp}
\ee
The power of 6 was chosen to provide a steep, but continuous
profile so that both the phase error and the oscillation error are visible. If the profile were too
steep, like that of a square wave, the oscillation error would have overwhelmed the calculation before
the phase error can be seen. The oscillation errors in all these symplectic algorithms 
are primarily due to the oscillation error in algorithm 1A. Algorithm 1B has a much smaller
oscillation error. Because all the algorithms are essentially norm-preserving, oscillation errors 
is inherent to any scheme which does not preserve the positivity of the solution\cite{hv}.

Higher order advection algorithms can again be constructed by the method of composition (\ref{t2n}).
Since $g_2$ is unitary, any product of $g_2$ is also unitary. 
Thus to preserve unitarity, one must use only a single product composition, rather than a multi-product 
expansion as in the diffusion case. (However, as noted in the last secrion, this violation of unitarity in MPE 
is small with increasing order. At sufficiently high order, this violation is beyond machine precision and is 
indistinguishable from a truly unitary algorithm\cite{bcr99}. For simplicity, we will only consider strictly 
unitary algorithms in this discussion.) For a single product composition,
the resulting phase angle is just a sum of $\phi_2$'s.  
The simplest fourth-order composition, the Forest-Ruth (FR) algorithm \cite{cre89,fr90,yos90} is given by  
\be
T_4^{F\!R}(\dt)=T_2(a_1\dt)T_2(a_0\dt)T_2(a_1\dt),
\la{fr}
\ee
with $a_1=1/(2-b)$, $a_0=-b/(2-b)$ and $b=2^{1/3}$. The coefficients $a_0$ and $a_1$ satisfy the consistency condition
$2a_1+a_0=1$ and the fourth-order condition $2a_1^3+a_0^3=0$. 
If we take $T_2(\dt)$ to be A2C, then the phase angle for $T_4^{F\!R}(\dt)$
is just (\ref{phi2cn}) with all $\eta^3$ terms removed, 
\be
\phi_4^{F\!R}=\eta\theta-\frac\eta{6}\theta^3+\left(\frac\eta{120}-5.29145\frac{3\eta^5}{240}\right)\theta^5+\cdots,
\la{cn4}
\ee
which is then correct to fourth-order in $\theta$. This is not true if we take $T_2(\dt)$ to be
the original algorithm A2. That fourth-order time-marching algorithm's phase angle
will still have a small error slope at $\theta=0$. One should therefore only uses A2C 
to compose higher order algorithms. 

The large numerical coefficient in (\ref{cn4}) is due to $2a_1^5+a_0^5=-5.29145$,
reflecting the fact that FR has a rather large residual error.
A better fourth-order algorithm advocated by Suzuki\cite{suzu90} (S4) at the expense of two more $T_2$ is
\be
T_4^S=T_2(a_1\dt)T_2(a_1\dt)T_2(a_0\dt)T_2(a_1\dt)T_2(a_1\dt)
\ee 
where now $a_1=1/(4-4^{1/3})$, $a_0=-4^{1/3}a_1$, and $4a_1^5+a_0^5=-0.074376$, which is nearly sixth-order.  

At the expense of two more $T_2$, one can achieve sixth-order 
via Yoshida's algorithm\cite{yos90} (Y6),
\be
T_6^Y=T_2(a_3\dt)T_2(a_2\dt)T_2(a_1\dt)T_2(a_0\dt)T_2(a_1\dt)T_2(a_2\dt)T_2(a_3\dt)
\ee
with coefficients 
\ba
a_1&=&-1.17767998417887\quad a_2= 0.235573213359357\nn\\
a_3&=& 0.784513610477560\quad{\rm and}\quad a_0=1-2(c_1+c_2+c_3). 
\la{yos}
\ea
For eighth and higher order algorithms, see Refs.\onlinecite{hairer02} and \onlinecite{mcl95}. 

The phase errors of these higher order algorithms at $\eta=0.7$ are shown at the right of Fig.\ref{stbad}. 
Since each algorithm applies $T_2$ (taken to be A2C) $n$ (=3,5,7) times, they are compared
to the phase error of running $T_2$ $n$ times at a reduce time step of $\dt/n$. Algorithms S4 and Y6
beat their target comparisons by orders of magnitude.  
There is a clear advantage in going to higher-order algorithms for solving the advection equation.   

\begin{figure}[hbt]
\includegraphics[width=0.7\linewidth]{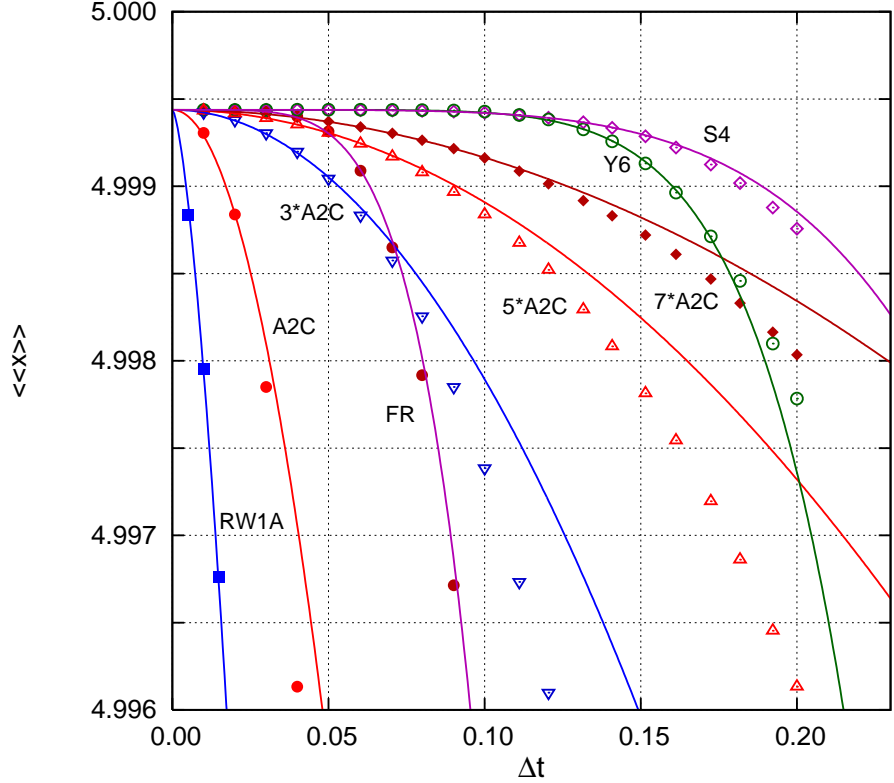}
\caption[]{\label{adcong} Comparing the convergence of various unconditionally stable symplectic 
advection algorithms. The solid lines are power laws of the form $a+b\dt^n$ seeking to verify 
the order of the algorithm. Higher order algorithms composed of $n$ second-order algorithms A2C
are compared to $n$*A2C at a reduce time-step size of $\dt/n$. FR, S4 and Y6 requires 
3, 5 and 7 runs of A2C respectively. In this computation of (\ref{xabs}), all algorithms met or 
exceeded their nominal order of convergence. See text for details.     
}
\end{figure}

In Fig.\ref{adcong}, the convergence of these higher order 
algorithms are compared. The range of the time steps used, 
$\dt=0.01-0.20$, corresponds to $\eta=0.4-8.0$. The same profile (\ref{intp}) is initially centered 
at $x=-5$ and propagated to $x=5$ at $v=1$. The time steps are chosen as 
$\dt=10/m$, so that $m$ iterations exactly give $t=10$. Since all algorithms satisfy (\ref{advec}) despite
the oscillation errors, we compute the expectation value
\be
\langle\langle x \rangle\rangle=\frac{\sum_{j=1}^N (j\dx)|u_j|}{\sum_{j=1}^N |u_j|}
\la{xabs}
\ee
with respect to the absolute value of the propagated profile. The phase errors for A2 and A2S
are known to be large from Fig.\ref{alg7}, and are not included in this comparison. The solid lines are
fitted power laws of the form $a+b\dt^n$. The first suprise is that
RW1A's result cannot be fitted with $n=1$. The fitted line is a fit with $n=3/2$.
The algorithm A2C can be well fitted with $n=2$. This is specially clear 
in the case where A2C is applied seven times at step size $\dt/7$. The fourth-order Forest-Ruth (FR)
and Suzuki (S4) algorithms can only be fitted with $n=5$, and the sixth-order Yoshida (Y6) algorithm
with $n=7$. They all converged to a value of $a=4.999438$ which is below the exact value of 5.
This is related to the grid size error. Halving the grid size to $\dx=0.0125$ gives
$a=4.999997$.

\section {Symplectic advection-diffusion algorithms}

The advection-diffusion equation,
\be
\frac{\partial u}{\partial t}=-v\frac{\partial u}{\partial x}+D\frac{\partial^2 u}{\partial x^2},
\la{addiff}
\ee
has the exact operator solution
\be	 
u(x,\dt)=\e^{ -v\dt\frac{\partial}{\partial x}+D\dt\frac{\partial^2 }{\partial x^2}}u(x,0).
\ee
If $v$ and $D$ are just constants, then 
since $[\frac{\partial}{\partial x},\frac{\partial^2 }{\partial x^2}]=0$, one has
\ba
u(x,\dt)
&=&\e^{ -v\dt\frac{\partial}{\partial x}}\e^{D\dt\frac{\partial^2 }{\partial x^2}}u(x,0)\nn\\
&=&\e^{ -v\dt\frac{\partial}{\partial x}}\widetilde u(x,\dt)\nn\\
&=&\widetilde u(x-v\dt,\dt),
\ea
where $\widetilde u(x,\dt)$ is the diffused solution. The complete solution is therefore
the exact diffused solution $\widetilde u(x,\dt)$ displaced by $v\dt$.

For periodic boundary condition, our matrices {\it also} commute, $[\A,\B]=0$, so that the discretized
version also holds,
\be
\bu(t+\dt)=\e^{\dt\A}\e^{\dt\B}\bu(t).
\la{abeq}
\ee
Thus arbitrary high order algorithms can be obtained by applying higher order advection 
and diffusion algorithms in turns from the previous sections. In the case of  
spatially dependent $D(x)$ or $v(x)$ where $[\A,\B]\ne 0$, one can do the second-order splitting
\be										 
\bu(t+\dt)=\e^{\frac12\dt\A}\e^{\dt\B}\e^{\frac12 \dt\A}\bu(t)
\la{abaeq}
\ee
and apply higher order MPE algorithms. However, for periodic boundary condition, this way of 
solving the advection-diffusion equation cannot conserve the norm. Consider the case of applying the 
advection algorithms RW1A, RW1B followed by any norm-conserving diffusion algorithm. 
From (\ref{nrw1a}), the change in the modified norm after $\dt$ would be 
\be
\tilde N_{\rm 1A}^\prime-\tilde N_{\rm 1A} =(\sqrt{1+\eta}-1)(u_1^\prime-u_1).
\la{chn1a}
\ee
For RW1A, $u_1$ is a discontinuous point {\it higher} than its adjacent neighbors $u_n$ and $u_2$. Consequently,
after the diffusion step, $u_1^\prime<u_1$ and there is a loss of normalization in (\ref{chn1a}). For
RW1B, $u_1$ is a discontinuous point {\it lower} than its adjacent neighbors $u_n$ and $u_2$. After
the diffusion step, $u_1^\prime>u_1$, and again results in a loss of normalization:
\be
\tilde N_{\rm 1B}^\prime-\tilde N_{\rm 1B} =(\sqrt{1-\eta}-1)(u_1^\prime-u_1).
\la{chn1b}
\ee
As will be shown, this loss is small for small $D$, but is irreversible and accumulative
after each orbit around the periodic box. 
For fixed boundary with $u_1=0$, there is no such norm-conserving problem. 
  
An alternative is to update the advection and diffuion steps simultaneously.
In this case, one might
decomposing $\A+\B$ into a sum of $2\times 2$ matrices as
done previously,
\be
\C_j=\A_j+\B_j =\frac{D}{\dx^2}\left(\baa{cccccc}
  \ddots&    &   & &  &   \\
  & -1 & 1 & &  &    \\
   &1 & -1& &  &    \\
   &     &  & & &\ddots     \\
		 \eaa\right)
  +\frac{v}{2\dx}\left(\baa{cccccc}
  \ddots&    &   & &  &   \\
  & 0 & -1 & &  &    \\
   &1 & 0& &  &    \\
   &     &  & & &\ddots     \\
		 \eaa\right)
\ee
resulting in
\be
\e^{\dt\C_j}=\left(\baa{cccccc}
  1&    &   & &  &   \\
  &\alpha &\lm & &  &    \\
   &\beta &\alpha& &  &    \\
   &     &  & & &1     \\
		 \eaa\right),
\quad\quad		 
\e^{\dt\C_N}=\left(\baa{cccccc}
  \alpha&    &   & &  &\lm   \\
  &1 &  & &  &    \\
   & & &1  &    \\
   \beta&     &  & & &\alpha    \\
		 \eaa\right),		
\ee
with
\be
\alpha=\e^{-r}\cosh\psi,
\quad\beta=\e^{-r}(r+\eta/2)\frac{\sinh\psi}{\psi},
\quad\lm=\e^{-r}(r-\eta/2)\frac{\sinh\psi}{\psi},
\la{dd}
\ee
and $\psi=\sqrt{r^2-(\eta/2)^2}$. The corresponding Saul'pev form of the 1A algorithm
is then
\be 
u_j^\prime 
=\beta u_{j-1}^\prime+\gamma u_j+\lm u_{j+1}
\la{add1a}
\ee
where $\gamma=\alpha^2-\beta\lm=\e^{-2r}$ remains the determinant of the updating matrix.
However, for the Saul'pev form (\ref{add1a}) to be norm-preserving, one must have
\be
\beta+\gamma+\lm=1.
\la{ncond}
\ee 
Surprisingly, this is grossly violated by (\ref{dd}) when both $r$ and $\eta$ are non-vanishing. 

As we have learned in the previous two sections, any such initial algorithm can be far from optimal.
Therefore, one may as well begin with an assumed updating matrix,
\ba
u_j^\prime&=& \alpha u_j+\lm u_{j+1}\nn\\
u_{j+1}^\prime&=&\beta u_j+\alpha u_{j+1}
\la{addup}
\ea
and determine its elements by enforcing norm-conserving condition (\ref{ncond}) and by matching the 
expansion coefficients of the exact amplification factor.
The sequential applications of this updating matrix yields Saul'pev-type algorithms 1A (\ref{add1a}) and 1B,
\be
u_j^\prime =\beta u_{j-1}+\gamma u_j+\lm u_{j+1}^\prime
\ee
The determinant $\gamma=\alpha^2-\beta\lm$ is to be regarded as fixing 
$\alpha$ as a function of $\gamma$ and $\beta$ via $\alpha=\sqrt{\gamma+\beta\lm}$.   
The resulting amplification factors are then
\ba
g_{\rm 1A}&=&\frac{\gamma+\lm\e^{i\theta}}{1-\beta\e^{-i\theta} }=\e^{-h_{\rm 1A}},\nn\\
g_{\rm 1B}&=&\frac{\gamma+\beta\e^{-i\theta}}{1-\lm\e^{i\theta} }=\e^{-h_{\rm 1B}}.
\ea
The norm condition (\ref{ncond}) fixes $\lm$ in terms of $\beta$ and $\gamma$. 
In terms of $\gamma$ and $\beta$ algorithms 1A and 1B have expansions,
\ba
h_{\rm 1A}&=&\frac{2\beta-(1-\gamma)}{1-\beta}i\theta+\frac{(1-\gamma)(\beta+\gamma)}{2(1-\beta)^2}\theta^2 +O(\theta^3)\nn\\
h_{\rm 1B}&=&\frac{2\beta-(1-\gamma)}{\gamma+\beta}i\theta+\frac{(1-\gamma)(1-\beta)}{2(\gamma+\beta)^2}\theta^2 +O(\theta^3)
\ea
Matching the first and second order coefficients of the exact exponent 
\ba
h_{ex}&=&i\eta\sin(\theta)+4r\sin(\theta/2)^2\nn\\
&=&i\eta\theta+r\theta^2 +O(\theta^3),
\la{addgex}
\ea
then completely determines, for 1A and 1B respectively,
\be
\beta= \frac{1-\gamma+\eta}{2+\eta}\quad \gamma=\frac{1-wr}{1+wr}\quad w=\frac2{2+\eta(3+\eta)},
\ee
\be
\beta= \frac{1-\gamma+\gamma\eta}{2-\eta}\quad \gamma=\frac{1-wr}{1+wr}\quad w=\frac2{2-\eta(3-\eta)}.
\ee
These are the generalized Roberts-Weiss algorithms for the advection-diffusion equation. 

\begin{figure}[hbt]
\includegraphics[width=0.7\linewidth]{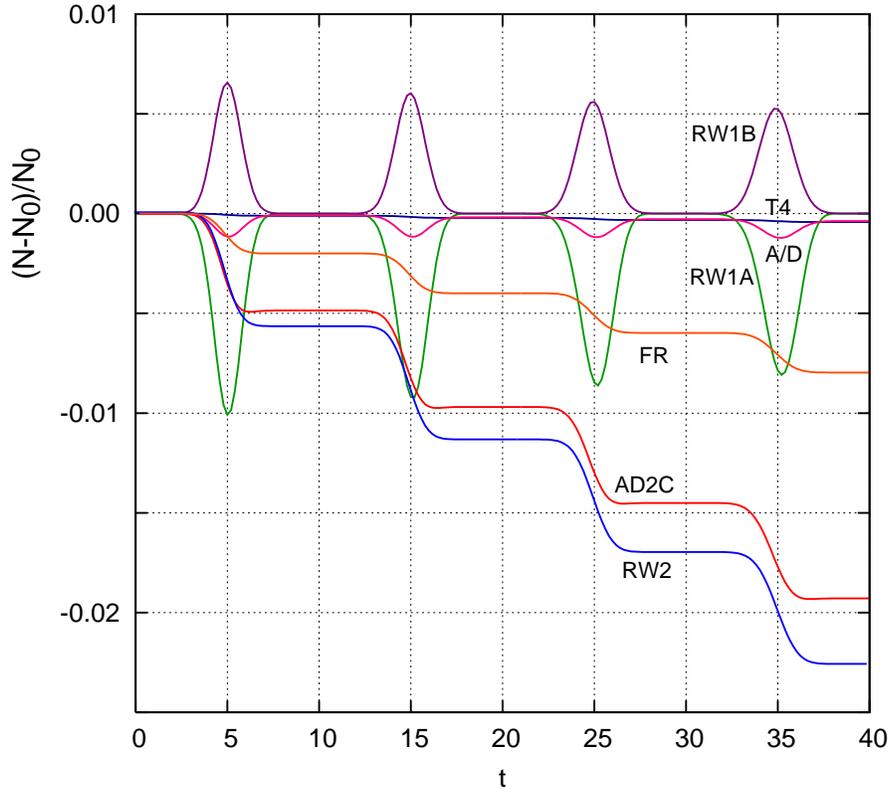}
\caption[]{\label{norm} 
The normalization error of various algorithms when propagating a Guassian profile in a periodic box of
[0,10] with $\dx=0.05$, $\dt=0.033$, $v=1$, $D=0.005$, $\eta=0.66$ and $r=0.066$. The profile is
initially centered at $x=5$. At $t=5,15,25,25$, the Gaussian peak is at the edge of the periodic box.
First-order advection-diffusion algorithms RW1A and RW1B conserve the norm periodically. All higher
than first-order algorithms suffer loss of normalization irreversibly, though very small for  
MPE algorithm T4 and A/D. The latter is applying the advection algorithm A2C
and the diffusion algorithm D2S sequentially.
}
\end{figure}

For any choice of $\beta$ and $\gamma$, the modified norm including the boundary effect now reads
\ba
\tilde N_{\rm 1A}&=&N+(\frac{\alpha}{1-\beta}-1)u_1\la{addn1a}\\
\tilde N_{\rm 1B}&=&N+(\frac{\alpha}{1-\lm}-1)u_1
\la{addn1b}
\ea
Remarkably, for the above generalized RW algorithms, one has
\be
\frac{\alpha}{1-\beta}=\sqrt{1+\eta}\quad{\rm and}\quad \frac{\alpha}{1-\lm}=\sqrt{1-\eta}.
\ee 
The modified norms (\ref{addn1a}) and (\ref{addn1b}) are therefore the same as 
the pure advection cases of (\ref{nrw1a}) and (\ref{nrw1b}). For these two first order
advection-diffusion algorithms, their norm-conservation in a periodic box 
will then be periodic, as in the pure advection case. This is shown in Fig.\ref{norm}.
However, as soon as one concatenate them into algorithm RW2, the loss of norm is irreversible.
This is because each algorithm will behave as a diffusion algorithm for the other. 
The norm-loss mechanism described earlier will then apply. 

The second-order algorithm's amplification factor is
\be
g_2=\left(\frac{\widetilde\gamma+\widetilde\lm\e^{i\theta}}{1-\widetilde\beta\e^{-i\theta} }\right)
\left(\frac{\widetilde\gamma+\widetilde\beta\e^{-i\theta}}{1-\widetilde\lm\e^{+i\theta} }\right)=\e^{-h_2},
\ee
where $\widetilde\gamma=\gamma(\dt/2)$, etc..	In terms of $\widetilde\gamma$ and $\widetilde\beta$, 
$h_2$ has the expansion,
\be
h_2=i\theta
\left(\frac{(1+\widetilde\gamma)(\widetilde\gamma-1+2\widetilde\beta)}
{(1-\widetilde\beta)(\widetilde\gamma+\widetilde\beta)}\right)+ O(\theta^2).
\ee
Matching this to the first order coefficient of the exact exponent (\ref{addgex}) 
determines
\be
 \widetilde\beta=\frac12(1- \widetilde\gamma)+\frac12(1+\widetilde\gamma)\widetilde s,
\ee
and
\be
 \widetilde\lm=\frac12(1- \widetilde\gamma)-\frac12(1+\widetilde\gamma)\widetilde s.
\ee
where $\widetilde s$ has been previously defined by (\ref{stilde}).
In terms of only $\widetilde\gamma$,
\be
h_2=i\eta\theta +2\frac{(1-\widetilde\gamma)}{(1+\widetilde\gamma)}\frac{(1+3\widetilde s^2)}{(1-\widetilde s^2)^2}\theta^2
+O(\theta^3),
\ee
and matching the second order coefficient in (\ref{addgex}) determines
\be
\widetilde\gamma=\frac{1-w r/2}{1+w r/2}\quad{\rm with}\quad w=(1-\widetilde s^2)^2/(1+3\widetilde s^2).
\ee
If $\eta=0$, $\widetilde s=0$, one recovers (\ref{saulcoef}), which is the second-order diffusion
algorithm D2S. If $r=0$, then $\gamma=1$ and one recovers the second-order advection algorithm A2C with   
$\widetilde\lm=-\widetilde\beta=-\widetilde s$. We shall refer to this second-order algorithm as AD2C. 
AD2C is an unconditionally stable algorithm which requires only half the effort of algorithm A/D,
which applies A2C and D2S sequentually. However, AD2C has a greater irreversible norm-error 
when applied to periodic boundary problems. This is shown in Fig.\ref{norm}.

\begin{figure}[hbt]
\includegraphics[width=0.7\linewidth]{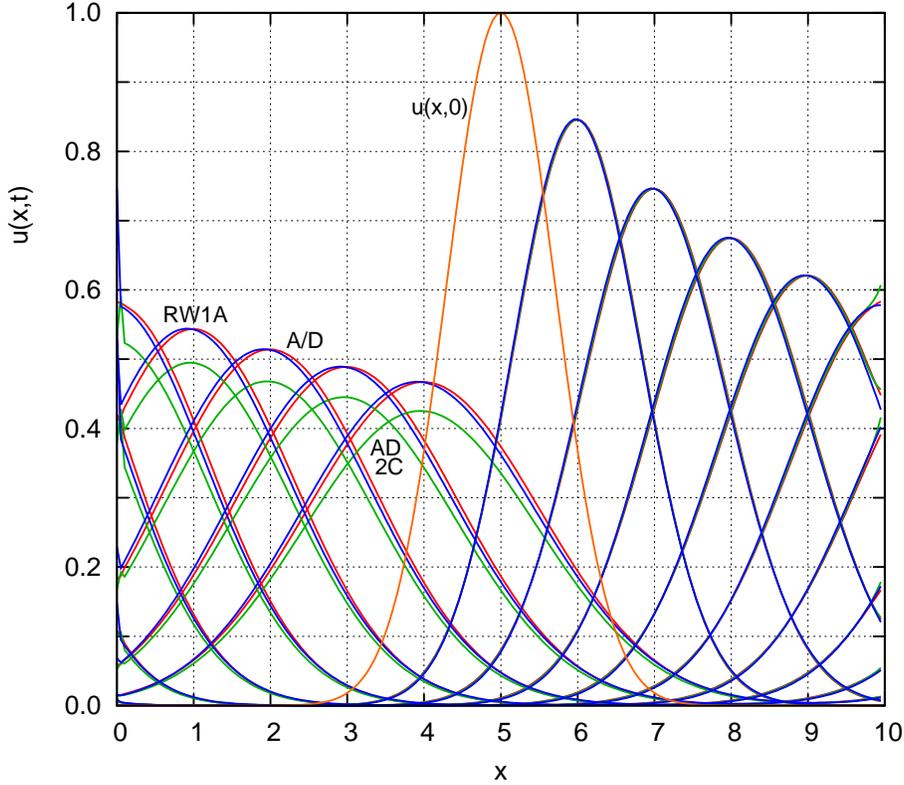}
\caption[]{\label{add} 
The propagation of a Guassian profile in a periodic box of
[0,10] with $\dx=0.05$, $\dt=0.033$, $v=1$, $D=0.1$, $\eta=0.66$ and $r=1.33$. The profile is
initially centered at $x=5$. All three profiles produced by algorithms RW1A, AD2C and A/D are in
essential agreement prior to the pulse peak hitting the right periodic edge. As the profiles reappear
from the left, the norm of AD2C is noticeably lower.
}
\end{figure}

This norm-error for periodic boundary condition can be greatly reduced by going to higher orders.   
Fig.\ref{norm} shows the result for the fourth-order MPE algorithm $T_4$, with AD2C as $T_2$. For
$r$ small, surprisingly, even the negative-coefficient algorithm FR is stable, but with an error 
comparable to second-order algorithms.

In Fig.\ref{add} we illustrate the effect of this norm-loss error at a large value of $r=1.33$. 
For clarity, only results from three representative algorithms are shown. Algorithms A/D and RW1A have
small or only periodic norm-losses and remained in agreement after the Gaussian peak has reappeared 
from the left. However, algorithm AD2C suffers an irreversible norm-loss and its peak is noticeably lower.

 \section {Concluding Summary}

In this work, we have shown that explicit symplectic finite-difference methods can be derived in the
same way as symplectic integrators by exponential splittings. The resulting sequential updating 
algorithms reproduce Saul'pev's unconditionally stable schemes, but is more general and 
can be applied to periodic boundary problems. In contrast to Saul'pev's original approach, where the 
algorithm is fixed by its derivation, symplectic algorithms can be systematically improved by matching 
the algorithm's amplification factor more closely to the amplification factor of the 
semi-discretized equation. One key contribution of this work is the recognition that, 
for finite difference schemes, their amplification factors should be compared, not to the 
continuum growth factor, but to the amplification factor of the semi-discretized equation. 
The exponent of this amplification factor then serve as  
the ``Hamiltonian" for developing symplectic finite-difference algorithms. By requiring the algorithm's 
modified ``Hamiltonian" to match the original ``Hamiltonian" to the leading order, one produces 
all known, non-pathological first-order Saul'pev schemes and many new second-order algorithms for 
solving the diffusion and the advection equation. As a consequence of this formal correspondence with 
symplectic integrators, existing methods of generating higher order integrators can be immediately
used to produce higher order finite-difference schemes.

The generalization to higher dimensions can be done by dimensional splitting, resulting in
unconditionally stable, alternate-direction-explicit methods. The generalization to 
non-constant diffusion and advection coefficients is a 
topic suitable for a future study. The coefficients must frozen in such a way that one
can recover Saul'pev's asymmetric schemes from their more basic sequentual updatings.

\begin{acknowledgments}
This work is supported in part by the Austrian FWF grant P21924 and
the Qatar National Research Fund (QNRF) National Priority Research Project (NPRP)
grant \# 5-674-1-114. I thank my colleague Eckhard Krotscheck and the Institute for Theoretical 
Physics at the Johannes Kepler Univeristy, Linz, Austria, for their wonderful hospitality 
during the summers of 2010-2012.
\end{acknowledgments}

\end{document}